\documentclass{ifacconf}

\usepackage{algorithmic}
\usepackage{algorithm}
\usepackage{natbib}
\usepackage{amsmath}
\usepackage{amssymb}
\usepackage{graphicx}

\newcommand{\PWL}[1][{}]{ (p,Q,\{(n_{q}^{#1}, A_{q}^{#1},a_{q}^{#1},C_{q}^{#1},c_{q}^{#1}, \mathcal{X}_{q,0}^{#1}) \mid q \in Q \})}

\newcommand{\PWLS}{\textbf{PWL}}
\newcommand{\SPWLS}{\textbf{PWLs}}
\newcommand{\SARS}{\textbf{SARS}}
\newcommand{\SAR}[1][{}]{(p,n,Q^{#1},\{ A^{#1}_{q,i} \mid q \in Q,i=1,\ldots,n \})}
\newcommand{\LSSJ}{\textbf{LSSJ}}
\newcommand{\SLSSJ}{\textbf{LSSJs}}
\newcommand{\SPAN}{\textrm{Span}}
\newcommand{\Rank}{\mathrm{rank}\mbox{ } }
\newcommand{\CLUST}{\mathbf{C}}

\newtheorem{Theorem}{Theorem}

\newtheorem{Proposition}{Proposition}
\newtheorem{Corollary}{Corollary}

\newtheorem{Definition}{Definition}
\newtheorem{Remark}{Remark}

\newtheorem{Example}{Example}

\newtheorem{Problem}{Problem}

\newcommand{\QNUM}{D}

\begin{document}
\begin{frontmatter}
\title{Identification of Piecewise Linear Models of Complex Dynamical Systems}

 \author{Ronald L. Westra}, 
 \author{Mih\'aly  Petreczky}
 \author{Ralf L.M. Peeters}

 \address{DKE, Maastricht University, The Netherlands, \{Westra,M.Petreczky,Ralf.Peeters\}@maastrichtuniversity.nl}

\begin{abstract}
 The paper address the realization and identification problem
 for a subclass of piecewise-affine hybrid systems. The paper
 provides necessary and sufficient conditions for 
 existence of a realization, a characterization of minimality, 
 and an identification algorithm for this
 subclass of hybrid systems. The considered system class and
 the identification problem are motivated by applications in
 systems biology.
\end{abstract}


\end{frontmatter}

\section{ Introduction}

In this paper we address the realization and identification problem
for a subclass of piecewise-affine hybrid systems.

\textbf{Contribution of the paper}
 We define the class of 
 \emph{piecewise-linear systems} (abbreviated by \PWLS).
 \PWLS{s} are a subclass of piecewise-affine hybrid systems.
 The continuous dynamics of \PWLS\ is determined by a finite
 collection of affine subsystems. However,
 in contrast to traditional
 piecewise-affine systems,
 we allow any change of the continuous state during a discrete-state transition, as long as the
new state belongs to the set of designated initial states of the
affine subsystem associated with the new discrete state.
In addition, we do not impose any specific mechanism
for triggering discrete-state transitions.

 We formulate the realization problem for this system class and
 partially solve it by providing necessary and sufficient conditions for
 existence of a realization. We also present conditions for minimality.
 We show that the outputs of any \PWLS\ can also be described by a switched AR 
 model.
 The main conclusion is that any system can be transformed to a minimal
 system \emph{with one discrete state} while preserving input-output behavior. 
This means
that without further restrictions, the identification problem for
such systems is not necessarily interesting. We discuss a number of
restrictions on the system structure which avoid this problem.

Note that the conclusion above is not valid for other classes of
hybrid systems. For hybrid systems from \cite{Petreczky09-TAC},
there examples of input-output maps which provenly cannot be
realized by systems with one discrete state.

 In addition, we present an identification algorithm for systems with
full observations. This algorithm is illustrated
by example of physical and biological relevance.


\textbf{Motivation}
 The motivation for studying realization theory for \PWLS{s} is that
 it provides the theoretical foundations for systems identification.
 The motivation for investigating identification of \PWLS{s} is the
 following.  
 \begin{itemize}
 \item
  First, several systems of interest can be modeled
 by \PWLS{s} and \PWLS{s} tend to be convenient for analysis. However
 the parameters of
 \PWLS\ models are often not directly available and hence they have to
 be estimated from measurements.
 \item 
    The problem of estimating the network dynamics of complex biological
    systems can be reduced to the identification problem for \PWLS{s}.
 \end{itemize}

\textbf{Identification of \PWLS{s} and biological networks}
 Below we will elaborate on the relationship between identification of
 \PWLS{s} and the estimation of dynamic interactions of complex dynamical
 systems.   

Numerous vital processes in nature involve complex signaling networks, 
varying from gene-protein interaction networks \cite{Westra2007} and 
complex communication in microbes \cite{Charu2006} to the synchronization 
in the heart of higher animals \cite{Heijman2009}. In most cases, 
these signaling networks do not have a fixed topology, but vary their 
structure according to their internal states and certain external conditions 
\cite{Adami2000}. This flexibility allows the organism to exhibit complex 
behavior and flexible responses to various environmental conditions.

The underlying biological problem is to describe these flexible dynamical 
networks, based 
on (partial) observations of their internal states. Examples of
potential observations are time series of 
gene expressions, protein densities, oxygen stress and sugar concentrations. 

\vspace{-5pt}

In almost all cases, there is no explicit mathematical model available that 
even remotely describes the network dynamics. This is caused by the sheer size 
of the problem; the many thousands of genes, RNAs, proteins and other types of 
molecules involved, the complexity and idiosyncrasy of their interactions, and 
the large amount of noise in biological systems. Basically, many agents 
involved in the network are unknown let alone their mechanism of interaction.

Mathematically, the processes above can be viewed as nonlinear dynamical
systems with unknown parameters. Such systems often have several equilibria,
and hence their behavior
can be seen as a mixture of certain affine systems, where each affine system is obtained by linearizing the original system around one of the  equilibria.
That is, such systems can be approximated by a \PWLS.
In addition to their simplicity, the advantage of \PWLS\ approximations is
that they neatly capture the interaction
among various state components around various equilibria. 
More precisely, interaction between the $i$th and $j$th state
components can be
viewed as the property that the $(i,j)$th entry of the Jacobian
of the system around equilibrium is non-zero.
Since we are often interested exactly in interactions rather than the
detailed dynamics, it makes sense to recast the problem 
of identifying such interactions into the problem of identifying 
\PWLS{s} approximations. While the paper
does not directly address identification of interactions, 
we believe that the presented results represent a step towards the solution of that problem.

\textbf{Related work}
  Identification of piecewise-linear (-affine) hybrid systems has been
  subject of intensive research,
  \cite{JLSVidal,MaVidal,Bako09-SYSID2,IdentSurvey,VidalAutomatica,PLClustering,Juloski1,Bako08-IJC,Roll:Automatica04,Fox:09}.
  Realization theory for hybrid systems was investigated in
  \cite{Weiland06,PaolettiTac,MP:Phd,Petreczky09-TAC,GrossmanHybAlg}. 
 The application of hybrid systems to modelling complex dynamical systems,
 originating from biology, and their identification
 has been a subject of intensive research
 \cite{PLGenRegStab,LitOverview,
TrecateBioIdent2,KouCin:2007:IFA_2865,CinEtal:2008:IFA_2976}.

The results of this paper are new, to the best of our knowledge.
What distinguishes the contribution of this paper from the existing
work is \textbf{(a)} the system class considered, in particular,
the emphasis on continuous-time systems, 
\textbf{(b)} the emphasis on realization theory and algorithms,
\textbf{(c)} the details of the identification algorithm.

\textbf{Outline of the paper}
\S \ref{sect:pwl} presents the formal definition of the system
class of interest.
\S \ref{sect:approx} discusses the relationship between nonlinear
systems with complex dynamics and piecewise-affine hybrid systems.
 \S \ref{sect:real} presents the results on realization
theory, and \S \ref{sect:ident} presents the identification algorithm
and the results of the numerical experiments.

\textbf{Notation}
  We use the standard notation. We denote by
 $T=[0, +\infty)$ the time-axis. We denote by $I_n$ the $n \times n$
 identity matrix. We denote the set of natural numbers including zero by
 $\mathbb{N}$.

\vspace{-10pt}

\section{Piecewise-Linear systems (PWL)}
\label{sect:pwl}
 The aim of the section is to define the class of piecewise-linear
 systems formally.
 \begin{Definition}[PWL]
  A \emph{piecewise-linear system} (abbreviated as PWL) is a 
  dynamical system determined by 
  \begin{equation}
  \label{phas:def}
   \Sigma \left\{\begin{split}
    & \dot x (t)=A_{q(t)}x(t)+a_{q(t)} \\
    &  y(t)=C_{q(t)}x(t)+c_{q(t)} \\
    &  x(t^{+}) \in \mathcal{X}_{q(t^{+}),0} \\
  \end{split}\right.
  \end{equation}
  Here $Q=\{1,\ldots,\QNUM\}$ is the finite set of discrete modes,
  $\mathbb{R}^{p}$ is the output space, $\mathbb{R}^{n_q}$ is the state-space
  of the system in mode $q \in Q$,
  For each $q, \in Q$
  $A_{q} \in \mathbb{R}^{n_q \times n_q}$, 
  $a_q \in \mathbb{R}^{n_q}$, 
  $C_q \in \mathbb{R}^{p \times n_q}$, and $c_q \in \mathbb{R}^{p}$
  are the parameters of the affine system in mode $q \in Q$,  
  $\mathcal{X}_{q,0} \subseteq \mathbb{R}^n$ -- is the set of 
  initial states of the affine system in mode $q \in Q$.
The \emph{state space} $\mathcal{H}_{\Sigma}$ of $\Sigma$ is
   \( \mathcal{H}_{\Sigma}=\bigcup_{q \in Q} \{q\} \times \mathbb{R}^{n_q} \)
  We call $\Sigma$ \emph{linear}, if $c_q=0$ and $a_q=0$ for all $q \in Q$,
  otherwise $\Sigma$ is called \emph{affine}.

   We will use the following short-hand notation
  \[ \Sigma=\PWL. \]
 \end{Definition}
 Informally,
 the evolution of $\Sigma$
 takes place as follows. 
 As long as the value of the discrete state $q(t)$ at time $t$ does not
 change,
 the continuous state and the
 continuous output an time $t$ change according to the
 affine system
 $\dot x(t) = A_{q(t)}x(t)+a_{q(t)}$ and $y(t)=C_{q(t)}x(t)+c_{q(t)}$.
 The discrete state can change at any time, however, we do not allow
 consecutive changes of discrete states immediately one after the other.
 If the discrete state changes to $q(t^{+})$ at time $t$, then 
 the new continuous state should satisfy $x(t^{+}) \in \mathcal{X}_{q(t^{+}),0}$,
 i.e. it should belong to the set of the designated initial states of
 the new discrete state $q(t^{+})$. 
 Note that we do not specify the mechanism which triggers the change of
 discrete state. 
 We also do not specify
 the initial discrete state, it is chosen by an unspecified mechanism.
 Hence, the \emph{description above allows for several state- and output-trajectories.}

 Note that the definition of \PWLS{s} presented above differs from the 
 one used in the literature.
 In the standard definition
 it is assumed that after a discrete-state transition, the new
 continuous state depends on the previous one.
 In contrast, here we only require that
 the new continuous state belongs to the designated set of
 initial states of the affine system associated with the new discrete state. This implies that \emph{the behavior of the system in discrete state
 $q$ depends only on the affine system associated with that discrete state and it does not depend on the discrete modes which were visited in the past.}
 Another consequence of the definition is that the model is not predictive.
 That is, the knowledge of the system parameters and the inputs does not
 uniquely determine the state- and output-trajectory of the system. 

 In order to define the evolution of \PWLS{s} formally, we
 need to introduce the following notation and terminology.
 \begin{Definition}[\cite{PieterHybTraject}]
  A time event sequence is a strictly monotone sequence
  $(t_{n})_{n=0}^{n^{*}}$ such that
  $n^{*} \in \mathbb{N} \cup \{+\infty\}$, $t_0=0$ and for all $0< n < n^{*}$,
  $0 \le t_{n} < t_{n+1}$. If $n^{*}=+\infty$ then we require that
   $\sup\{ t_n \mid n \in \mathbb{N}\} = +\infty$.
   If $n^{*} < +\infty$, then by convention $t_{n^{*}+1}=+\infty$.
 \end{Definition}
  The role of time event sequences is to formalize the time instances
  at which discrete events occur. The restrictions formulated in
  the definition imply that the set of switching times does not have
  an accumulation point, i.e. no Zeno-behavior can take place.
 \begin{Definition}[State-trajectory]
 A state-trajectory of $\Sigma$ is a map
 $\xi:T \rightarrow \mathcal{H}_{\Sigma}$ such that
 there exists a time event sequence
 $(t_i)_{i=0}^{n^{*}}$ and a sequence of discrete modes 
 $(q_i \in Q)_{i=0}^{n^{*}}$ such that
 for all $0 \le i \le n^{*}$, $i \in \mathbb{N}$, it holds that
 for all $s \in [t_{i},t_{i+1})$, $\xi(s)=(q_i,x(s-t_i))$ and
 \[ \dot x(t) = A_{q_i}x(t) + a_{q_i} \mbox{ and } x(0) \in \mathcal{X}_{q_i,0}. \] 
 The time event sequence $(t_n)_{n=0}^{n^{*}}$ is called the 
 sequence of \emph{switching times} of the state-trajectory $\xi$.

 We \emph{denote by $BS(\Sigma)$ the set of all state-trajectories of $\Sigma$}. 
 \end{Definition}
 \begin{Definition}[Output-trajectory]
  An output-trajectory of $\Sigma$ is a map $y:T \rightarrow \mathbb{R}^{p}$
  such that the following holds.
  There exists a state-trajectory $\xi$ of $\Sigma$, such that
  $y(t)=\upsilon_{\Sigma}(x(t))$ for all $t \in T$. Here,
  $\upsilon_{\Sigma}$ is the \emph{readout-map} of $\Sigma$, defined as
  \[ \upsilon_{\Sigma}: \mathcal{H}_{\Sigma} \ni (q,x) \mapsto C_qx \in \mathbb{R}^{p}. \]
  We \emph{denote by $B(\Sigma)$ the set of all output-trajectories of $\Sigma$}.
 \end{Definition}

 In the sequel,
unless stated otherwise, \emph{$f$ denotes a function $f:T \rightarrow \mathbb{R}^{p}$}.
The definition above implies that the external behavior of a \PWLS\ is
exactly a function of this type.
\begin{Definition}[Realization]
 The function $f$ is said to be realized
 by \PWLS\ $\Sigma$, if $if$ is an output-trajectory of $\Sigma$, i.e.
 if $f \in B(\Sigma)$. In this case $\Sigma$ is
 called a realization of $f$.
\end{Definition}
\begin{Definition}
 The dimension of a  \PWLS\ $\Sigma$ is a tuple $(|Q|,n)$, 
 where $n=\sum_{q \in Q} n_q$ if $\Sigma$ is linear, and
 $n=\sum_{q \in Q} (n_{q}+1)$, if $\Sigma$ is affine.
\end{Definition}
 That is, 
 the first element of $\dim \Sigma$ 
 is the number of discrete states, the second
 element is the number of continuous state components. The additional
 dimension in the case of affine \PWLS{s} stems from the need to
 store the vectors $a_q$, $c_q$, $q \in Q$.
We use the following partial order relation on
$\mathbb{N} \times \mathbb{N}$. We say that
$(p,q) \in \mathbb{N}$ is smaller than or equal  $(r,s) \in \mathbb{N}$,
denoted by $(p,q) \le (r,s)$,
 if $p \le r$ and $q \le s$.
Note that the order relation $\le$ in $\mathbb{N} \times \mathbb{N}$ is indeed a partial order, it is not possible to compare all
elements of $\mathbb{N} \times \mathbb{N}$.
\begin{Definition}[Minimality]
 A \PWLS\ $\Sigma$ is said to be a minimal realization of $f$, 
  if for any \PWLS\ $\hat{\Sigma}$ which is a realization of $f$,
 $\dim \Sigma \le \dim \hat{\Sigma}$.
\end{Definition}
The definition above implies that if $\Sigma$ is minimal, then it
is possible to compare its dimension to the dimension of any
other realization of $f$. Since we work with a partial order on dimensions,the existence of a minimal system is not at all obvious.
\begin{Problem}[Realization problem for \PWLS{s}]
 Find conditions
 for existence of a \PWLS\ realization of $f$. Characterize minimal
 \PWLS{s}\ realizations of $f$. Find algorithms for computing a 
 \PWLS{s}\ realization from finite data.
\end{Problem}
 As it was indicated before, our motivation for
 studying the realization problem for \PWLS{s} is to lay the theoretical
 foundations for identification of \PWLS{s}. 
 We present below the formulation of the identification problem.
\begin{Problem}[Identification problem]
 Assume that the value of $f$
 and its derivatives up to order $r$, $r \ge 0$
 is measured at time instances $t_1 < \ldots < t_k$. Based on the
 (possibly noisy) data
 $\{f^{(l)}(t_i)\}_{i=1\ldots,k, l=1,\ldots,r}$, find 
 a \PWLS\ realization of $f$.
\end{Problem}

\section{Approximation by \PWLS{s}}
\label{sect:approx}

Our main motivation for studying identification of
\PWLS{s} is the following. The problem of identifying 
interaction networks of complex dynamical systems can be
reduced to the identification problem of \PWLS{s}.
Indeed, consider a partially observed non-linear system of the form
\begin{equation}
\label{dyn:sys}
 \dot x(t) = f(x(t)) \mbox{ and } y(t)=h(x(t)) 
\end{equation}
 where $f:\mathbb{R}^{N} \rightarrow \mathbb{R}^{N}$ and
 $h:\mathbb{R}^N \rightarrow \mathbb{R}^p$ are sufficiently
 smooth functions.
We are interested in the dynamics of the network structure of \eqref{dyn:sys}.
By \emph{the network structure of \eqref{dyn:sys} at time $t$} 
we mean the directed graph 
with nodes numbered by $1,\ldots,N$ and where an edge goes from
node $i$ to node $j$ if the $i$th component of $x(t)$ influences the
$j$th component of $\dot x(t)$. As the values of $x(t)$ changes, so
does $f(x(t))$. Hence the network structure described above depends
on time.

The basic question of interest is how to derive the dynamics of
the interaction network based on the observed output $y$.
One obvious approach is to estimate $f$ and $h$ from $y$, but
this is feasible only if special restrictions are put on $f$ and $h$.
Notice, however, that if we look at local linearizations of $f$ around an
equilibrium and the 
interactions are strong enough, then
$x_i(t)$ influences $\dot x_j(t)$, if the 
$(i,j)$th entry of the Jacobian of $f$ is non-zero. 

Prompted by the considerations above we propose the following approach
to modelling and identification of the network structure of systems of
the form \eqref{dyn:sys}. We assume that \eqref{dyn:sys} has finitely
many equilibrium points.  We assume that 
under influence of noise the system jumps between neighborhoods
of these equilibrium points. We assume that 
the transition from the neighborhood of
one equilibrium point to the neighborhood of another equilibrium point
takes place very fast.
Hence, the behavior of the system during the transition can be ignored.
We associate a
\PWLS\ $\Sigma$ of the form \eqref{phas:def}
with \eqref{dyn:sys} as follows. 
Let $\QNUM$ be the number of equilibrium points. For each
$q \in Q$, let $e_q$ be the corresponding equilibrium point and
define
\[ 
  \begin{split}
   & A_q=Df(e_q) \mbox{ and } a_q=-A_{q}e_q  \\
   & C_q=Dh(e_q) \mbox{ and } c_q= h(e_q) 
  \end{split}
 \]
 Here $Df$ and $Dh$ denote the Jacobian of $f$ and $h$ respectively.
 Let $\mathcal{X}_{q,0}$ be a suitably chosen neighborhood of
 $e_q$. It then follows that if the state $x(t)$ of \eqref{dyn:sys}
 is close enough to $e_q$, we can approximate \eqref{dyn:sys} as follows
 \[ 
   \begin{split}
    & \dot x(t) =f(x(t)) \approx A_{q}x(t) + a_q \\
    &   y(t) = h(x(t)) \approx C_{q}x(t) + c_q 
   \end{split}
 \]
  As it was noted above, the relevant information regarding the
  interactions of state variables is already included in $A_q$.
  Hence, \emph{for the purposes of network reconstructions, identification
  of the \PWLS{s} approximation is sufficient}.

\section{Realization theory of \PWLS{s}}
\label{sect:real}
 Below we present some basic results on realization theory of \PWLS{s}.
 First we present a characterization of minimality, after that we 
 present conditions for existence of a realization and a realization
 algorithm. Throughout the section, \emph{$f$ denotes a function $f:T \rightarrow \mathbb{R}^{p}$}.

 \subsection{Minimality}
 In order to present our first result, we introduce the notion of a
 \emph{linear system with state-jumps}.
\begin{Definition}
 A linear system with state jumps, abbreviated as (\LSSJ) is a 
 linear \PWLS\ $\Sigma$ of the form \eqref{phas:def} with one
 discrete state, i.e. $\QNUM=1$. 
 We will identify the \LSSJ\ $\Sigma$ with the collection of data
 $(n,C,A,\mathcal{X}_0)$, where $n=n_1,C=C_1,A=A_1$ and
 $\mathcal{X}_{1,0}=\mathcal{X}_0$.
\end{Definition}
 The reason we call a \PWLS\ with one discrete state a 
 linear system with state jumps is that the system behaves as a 
 linear system, with the exception that its state occasionally jumps
 back to one of the initial states. For \LSSJ{s}, we can easily
 define the concepts of observability and span-reachability.
 \begin{Definition}
  A \LSSJ\ $\Sigma=(n,C,A,\mathcal{X}_0)$ is called \emph{observable}, 
  if $(C,A)$ is an observable pair, and $\Sigma$ is called
  \emph{span-reachable}, if $\SPAN\{A^{k}x_0 \mid k=0,\ldots,n-1,x_0 \in \mathcal{X}_{0}\}=n$.
 \end{Definition}
 \begin{Remark}
 \label{LSSJ:min}
  Note that  $Z=\SPAN\{A^{k}x_0 \mid k=0,\ldots,n-1,x_0 \in \mathcal{X}_{0}\}$
  is an $A$ invariant subspace containing $\mathcal{X}_0$, 
  hence by restricting $A$ and $C$ to $Z$ we can transform $\Sigma$
  to a span-reachable \LSSJ\ $\Sigma_r$ while preserving output
  trajectories, i.e. $B(\Sigma_r)=B(\Sigma)$.  By applying linear
  observability reduction, we can transform the \LSSJ\ $\Sigma_r$ to
  a \LSSJ\ $\Sigma_m$ which is span-reachable, observable and
  has the same output trajectories as $\Sigma$, i.e.
  $B(\Sigma_m)=B(\Sigma)$.
 \end{Remark}
 \begin{Theorem}[Minimality]
 \label{phas:min}
  Assume that $f$ admits a realization by a \PWLS.
  Then $f$ admits a minimal \PWLS\ realization $\Sigma$ such that
  $\Sigma$ is a span-reachable and observable \LSSJ. A \LSSJ\ realization
  of $f$ is minimal if and only if it is span-reachable and observable.
 \end{Theorem}
 The theorem above says that in general, the external behavior of any
 \PWLS{s} can be represented by a linear system with several initial
 states. That is, without further restrictions, the realization and
 identification problems for \PWLS{s} are equivalent to that of 
 \SLSSJ.
 
 The proof of Theorem \ref{phas:min} relies on the following 
 transformations.
   \begin{Definition}[\PWLS\ to linear \PWLS]
   \label{pahs2linpahs}
    Define the 
    \emph{linear} \SPWLS\ $L(\Sigma)$ associated with an affine $\Sigma$ as follows.
    \[ L(\Sigma)=\PWL[L] \]
     where $n_{q}^{L}=n_q+1$, $c_q^{L}=0$, $a_q^L=0$ and
    \[
      \begin{split}
       & A^{L}_q=\begin{bmatrix} A_q & a_q \\ 0 & 0 \end{bmatrix} 
         \mbox{,  } 
          C^{L}_q=\begin{bmatrix} C_q & c_q \end{bmatrix}   \\
        & \mathcal{X}_{q,0}^{L}=\{ (x^T,1)^T \in \mathbb{R}^{n_q+1} \mid x \in \mathcal{X}_{q,0}  \}.
       \end{split}
      \]
   \end{Definition}  
   \begin{Proposition}
    The output trajectories of $\Sigma$ and $L(\Sigma)$
    coincide, i.e. $B(\Sigma)=B(L(\Sigma))$, and
    $\dim \Sigma = \dim L(\Sigma)$.
   \end{Proposition}
   \begin{pf}[Sketch]
    A state trajectory $\xi$ of $L(\Sigma)$ is always a map
    of the form $\xi(t)=(q(t),(x(t)^T,1)^T)$, $t \in T$ such that 
    $\hat{\xi}(t)=(q(t),x(t))$, $t \in T$ is a state-trajectory of
    $\Sigma$. Since 
     $\upsilon_{L(\Sigma)}((q,(x^T,1)^T))=C_{q}x+c_q=\upsilon_{\Sigma}((q,x))$, the statement of the proposition follows.
   \end{pf}
  \begin{Definition}[Linear \PWLS\ to \LSSJ]
  Let $\Sigma$ be a linear \PWLS\ of the form \eqref{phas:def}.
  Define the \LSSJ\ $LS(\Sigma)$ associated with $\Sigma$ as follows.
  \[ LS(\Sigma)=(n,C,A,\mathcal{X}_0), \]
   where $n=\sum_{q \in Q} n_q$ and 
   \[
   \begin{split}
     & A = \begin{bmatrix} A_1 & 0 & \cdots & 0  \\
                       0   & A_2 & \cdots & 0  \\
                       \vdots & \vdots & \cdots & \vdots   \\
                        0     &  0     & \cdots & A_{\QNUM} \\
                      \end{bmatrix} \mbox{ and }
     C  =  \begin{bmatrix} C_1 & \cdots & C_{\QNUM} \end{bmatrix} \\
     & \mathcal{X}_0  =  \bigcup_{q \in Q} \hat{\mathcal{X}}_{q,0} \\
     & \hat{\mathcal{X}}_{q,0} =
      \{ (\underbrace{0,\ldots,0}_{q-1-times}, x, 0,\ldots, 0)^T 
         \mid x \in \mathcal{X}_{q,0} \}, \forall q \in Q \\
  \end{split}
  \]
  \end{Definition}
  \begin{Proposition}
   The output trajectories of $\Sigma$ and $LS(\Sigma)$ are the same,
   i.e. $B(\Sigma)=B(LS(\Sigma))$, moreover, 
   $\dim LS(\Sigma)=(1,n)$ where $(|Q|,n)=\dim \Sigma$.
  \end{Proposition}
  \begin{pf}[Sketch]
   For each $q \in Q$, let $\mathcal{X}_q$ be the subset of
   $\mathbb{R}^{n}$ of the form $z=(\underbrace{0,0,\ldots,0}_{q-1-times},x^T,0,\ldots,0)^T$, $x \in \mathbb{R}^{n_q}$.
   For each $z \in \mathcal{X}_q$ of the above form, let
   $\Pi_q(z)=x \in \mathbb{R}^{n_q}$.
   It then follows that $\xi(t)=(1,x(t))$ is a state-trajectory of
   $LS(\Sigma)$ if and only if
   $\hat{\xi}(t)=(q(t),\Pi_{q(t)}(x(t)))$ is a state-trajectory of
   $\Sigma$, where $q(t)=q$ for some
   $q \in Q$ if and only if $x(t) \in \mathcal{X}_q$.
   Since $\upsilon_{LS(\Sigma)}((1,x))=C_q\Pi_q(x)=\upsilon_{\Sigma}((q,\Pi_q(x))$ for all $x \in \mathcal{X}_q$, the statement of the proposition follows.
  \end{pf}
  
  Using the two transformations defined above, and Remark \ref{LSSJ:min},
  we can present the proof of Theorem \ref{phas:min}.
 \begin{pf}[Proof of Theorem \ref{phas:min}]
  Assume that $\Sigma$ is a \PWLS\ realization of $f$.
  If $\Sigma$ is affine, then replace $\Sigma$ with the associated
  linear \PWLS\ $L(\Sigma)$, which is also a realization of $f$.
  Hence, we can assume that $\Sigma$ is already
  a linear \PWLS. Construct then the \LSSJ\  $LS(\Sigma)$ associated
  with $\Sigma$ and apply Remark \ref{LSSJ:min} to $LS(\Sigma)$ to
  obtain a span-reachable and observable \LSSJ\ $\Sigma_m$ such that
  $B(\Sigma_m)=B(LS(\Sigma))=B(\Sigma)$. It then follows that
  $\dim \Sigma_m \le \dim LS(\Sigma) \le \dim \Sigma$.

  Hence, since $\Sigma$ was arbitrary, it is enough to look for
  minimal realizations among span-reachable and observable \LSSJ{s}.
  Notice that for any two \LSSJ{s}\ $\Sigma_1$ and $\Sigma_2$, 
  either $\dim \Sigma_1 \le \dim \Sigma_2$ or 
  $\dim \Sigma_2 \le \dim \Sigma_1$, hence, among all the
  possible \LSSJ{s}\ realizations of $f$, there must exist a minimal one.
 \end{pf}

 \subsection{Existence of a realization}
 The conditions for existence of
 a realization will be formulated using the rank of the Hankel-matrix of $f$.
  In order to define the Hankel-matrix of $f$, we have to define
 the notion of Markov-parameters. To that end, we need the notion
 of \emph{piecewise-analytic functions}.
 \begin{Definition}[Piecewise-analytic]
  The map $f$ is called 
  \emph{piecewise-analytic}, if
  there exist a finite or infinite
 number of time instances, $t_i \in T$, $t_{i} < t_{i+1}$, 
 $i \le N_f$, $i \in \mathbb{N}$ 
 for some $N_f \in \mathbb{N} \cup \{+\infty\}$, such that the
 following holds. If $N_f < +\infty$, let $t_{N_f+1}=+\infty$. Then we
 require that $t_0=0$ and $\bigcup_{i=0}^{N_f} [t_i,t_{i+1}) = T$ and
 for each $i \in \mathbb{N}$, $i \le N_f$, 
 $f$ is analytic on $[t_i,t_{i+1})$, but $f$ is not analytic on any neighborhood of $t_i$ in $\mathbb{R}$. 
We call
 the points $\{t_i\}_{i=0}^{N_f}$ the points of \emph{non-analyticity}.
 We define the set
 \[ I_f = \{ i \in \mathbb{N} \mid i \le N_f \} \]
 of indices of points of non-analyticity.
 \end{Definition}
  The intuition behind the definition is as follows.
  If $f$ has a realization by a \PWLS\ $\Sigma$, then the only points where $f$
  is not analytic are the points where the corresponding state-trajectory
  of $\Sigma$ switches from one discrete mode to another. In fact, one can
  show that there always exists a state-trajectory of $\Sigma$ which yields
  $f$ as output trajectory and which switches only at time instances at which
  $f$ is not analytic.  Hence, if $f$ has a realization by a \PWLS, then
  $f$ is piecewise-analytic and the points of non-analyticity tell us the
  switching times of a \PWLS\ realization of $f$.

  In the sequel, \emph{$f$ is assumed to be piecewise-analytic}.
 \begin{Definition}[Markov-parameters]
  Assume that $f$ is piecewise-analytic and let
  $\{t_i\}_{i=0}^{N_f}$ be the points of non-analyticity of $f$.
 For each $i \in I_f$, define the $i$th
 \emph{Markov-parameter} $\mathbf{M}^f_i$ of $f$ as a sequence
 \( \mathbf{M}^{f}_i:\mathbb{N} \rightarrow \mathbb{R}^{p} \)
  \[ \forall k \in \mathbb{N}: \mathbf{M}^{f}_i(k)=\frac{d^{k}}{dt^k} f(t_i+s)|_{s=0}. \]
 \end{Definition}
  It is easy to see that the collection of
  Markov-parameters $\{\mathbf{M}_i^f\}_{i=0}^{N_f}$ determines
  the map $f$ uniquely. We use the Markov-parameters to
  define the Hankel-matrix of $f$.
  \begin{Definition}[Hankel-matrix]
  We define the \emph{Hankel-matrix} $H_f$ of $f$ as the infinite matrix,
  rows of which are indexed by $\mathbb{N} \times \{1,\ldots,p\}$, and 
  columns of which
  are indexed by $\mathbb{N} \times I_f$.
  The entry of $H_f$ indexed by row index $(i,r)$ and by column index
  $(j,l)$ equals 
  \[ [H_f]_{(i,r),(j,l)}=(\mathbf{M}_{l}^f(i+j))_{r} \]
   where
  $(\mathbf{M}_{l}^{f}(i+j))_{r}$ denotes the $r$th entry of 
  Markov-parameter $\mathbf{M}^{f}_{l}(i+j)$. The rank of $H_f$,
  denoted by $\Rank H_f$, is the dimension of the linear space
  spanned by the columns of $H_f$.
  \end{Definition}
  \begin{Theorem}[Existence of a \PWLS\ realization]
  \label{theo:ex}
   The map $f$ can be realized by a \PWLS\ if and only if it is
   piecewise-analytic and $\Rank H_{f} < +\infty$. Moreover, a
   minimal \LSSJ\ realization of $f$ of dimension $(1,\Rank H_f)$
   can be constructed from $H_f$ .
  \end{Theorem}
  \begin{pf}[Sketch]
   \textbf{only if\ } 
   Assume that $f$ can be realized by a \PWLS. Then it can be
   realized by a \LSSJ\ $\Sigma=(n,C,A,\mathcal{X}_0)$. Then
   $M_i^f(k)=CA^{k}x_{i}$ $\forall i \in I_f, k \in \mathbb{N}$
    for some $x_i \in \mathcal{X}_0$, i.e. 
    $\{M_i^f(k)\}_{k=0}^{\infty}$ are the Markov-parameters
   of the linear system $(C,A)$ from some initial condition. 
   It then follows from linear systems 
   theory that $\Rank H_f \le n < +\infty$.

  \textbf{if \ } 
   Assume that $n=\Rank H_f$ and
   fix a basis in the linear span of the columns of $H_f$.
   Define $C$ as the matrix in this basis of the 
   linear map which maps a column to
   the vector 
   formed by the rows of that column indexed
   by $(0,1),\ldots, (0,p)$, in this order.
   Define the matrix $A$ as the matrix in this basis of the linear map
   which maps the column indexed by $(j,l)$ to the column indexed
   by $(j+1,l)$. Finally, let $\mathcal{X}_0$ be the set of
   coordinate vectors of the columns of $H_f$ indexed  by $(0,l)$, 
   $l \in I_f$. Then $\Sigma_f=(n,C,A,\mathcal{X}_0)$ is a \LSSJ\ 
   realization of $f$. If $\Sigma$ is another \LSSJ\ realization of $f$,
   then from the \textbf{only if} part it follows that
   $\Rank H_{f}=\dim \Sigma_f \le \dim \Sigma$, i.e. $\Sigma_f$ is a minimal
   realization of $f$.
  \end{pf}
  
  In Algorithm \ref{alg0} we present a Kalman-Ho-like realization algorithm for
  \LSSJ{s}\ and hence \PWLS{s}. To this end, we need the following definition.
  \begin{Definition}[Finite Hankel sub-matrix]
  \label{def:part}
  Fix integers \\
  $R,L,M> 0$. 
  Define the set
  \[ I_f^R = \{ i \in I_f \mid i \le R\}. \]
  Define the submatrix $H_{f,L,M,R}$ of $H_f$ as the matrix
  which is formed by the intersection of
  the rows of $H_f$ indexed by the elements of
  $I_L=\{0,\ldots,L\} \times \{1,\ldots,p\}$ and the columns of $H_f$ indexed
  by the elements of $J_{M,R}=\{0,\ldots,M\} \times I_f^R$. 
 \end{Definition} 
   \begin{algorithm}
   \caption{
    \newline
    \textbf{Inputs:} Hankel-matrix $H_{f,L,M+1,R}$.  
    \newline
    \textbf{Output:} \LSSJ\ $\Sigma_{L,M,R}$.
    }
    \label{alg0} 
    \begin{algorithmic}[1]
    \STATE
        Compute the decomposition
        $H_{f,L,M+1,R} = \mathbf{O}\mathbf{R}$ such that 
        $\mathbf{O} \in \mathbb{R}^{I_{L} \times n}$ and
        $\mathbf{R} \in \mathbb{R}^{n \times J_{M+1,R}}$ and 
        $\Rank \mathbf{R} =\Rank \mathbf{O} = n$.
   \STATE
     Define $\bar{\mathbf{R}} \in \mathbb{R}^{n \times J_M}$ as the
     matrix formed by the columns of $\mathbf{R}$ indexed by elements of
     $J_{M,R}$. Let $\hat{\mathbf{R}}\in \mathbb{R}^{n \times J_{M,R}}$  
     the matrix such that 
     the column of $\hat{\mathbf{R}}$ indexed by $(j,l)$ equals 
     the column of $\mathbf{R}$ indexed by $(j+1,l)$.

    \STATE
       
      Define $\Sigma_{L,M,R}=(n,C,A,\mathcal{X}_0)$ as follows.
      \begin{itemize}
      \item
           The $i$th row of $C$ equals the row of $\mathbf{O}$ indexed by
           $(i,0)$, $i=1,\ldots,p$.
      \item
          The matrix $A$ is the solution of the equation
          \[ \hat{\mathbf{R}}=A\bar{\mathbf{R}} \]
      \item
          The set $\mathcal{X}_0$ consists of the columns of $\mathbf{R}$ indexed
          by indices of the form $(0,l)$, $l \in I_f^R$.
      \end{itemize}
    \end{algorithmic}
   \end{algorithm}
  \begin{Theorem}[Correctness of Algorithm \ref{alg0}]
  \label{alg0:theo}
   If $\Rank H_{f} \le n$, and
   $|\{M_{i}^{f} \mid i \in I_f\}| < +\infty$, 
   then for some $R > 0$, 
   the \LSSJ\ returned by Algorithm \ref{alg0} is a minimal realization of $f$.
   The conditions above hold if
   $f$ has a realization by a \LSSJ\ of dimension at most $n$ and with a
   finite set of initial states.
  \end{Theorem}

\subsection{Ill-posedness of the realization problem}
\label{sect:ill-pos}
 The results on realization theory of \PWLS{s} indicate that the
 realization and identification problems for \PWLS{s} are in general
 ill-posed in the following sense. Realizability by a \PWLS{s} is
 equivalent to realizability by a minimal
 \PWLS{s} with one single discrete state.
 Hence, there seems to be
 no intrinsic way to choose discrete states based on the data.
 Since the original intention was use the combination of 
 the dynamics in various discrete modes to explain complex dynamics, the
 above conclusion is an unpleasant one.

 Note that the main problem is that one can increase the dimension of
 the continuous state-space to encode discrete states. 
  One way to remedy this is to
  place restriction on the dimension. To this end, we introduce
  the notion of $K-N$ realization.
  \begin{Definition}[$K-N$]
  A \PWLS\ $\Sigma$ is said to be a $K-N$ \PWLS\
  if $\Sigma$ is of the form \eqref{phas:def} and $|Q| \le K$ and
  for all $q \in Q$, $n_q \le N$. The map $f$ is said to admit a
  $K-N$ realization, if there exists a $K-N$ \PWLS\ which is a 
  realization of $f$.
  \end{Definition}
   The motivation for using $K-N$ realizations is that
   by choosing an appropriate $N$, we can make sure that
   the growth in the number of continuous states cannot be used
   to replace discrete states.
 Below we present conditions for existence of a $K-N$ realization of
 $f$. To this end,
 we have to introduce the following concept.
 \begin{Definition}[Hankel-matrices \& partition] 
 Consider the partitioning $\CLUST=(\CLUST_q)_{q=1}^{K}$ of the set of $I_f$.
 For each $q=1,\ldots, K$, let
 $H_{f,\CLUST}^{q}$ be the sub-matrix of Hankel-matrix $H_f$ which is
 formed by the columns of $H_f$ indexed
 by indices of the form $(j,l)$, $j \in \mathbb{N}$ and $l \in \CLUST_q$.
 \end{Definition}
 \begin{Theorem}[Existence of $K-N$ realizations]
 \label{theo:pos1}
  The map $f$ has a 
  realization by a linear $K-N$ \PWLS\
  if and only if there exists a partitioning 
  $\CLUST=(\CLUST_i)_{i=1}^{\QNUM}$ of $I$ with $\QNUM \le K$
  such that for all $i=1,\ldots,\QNUM$,
  $\Rank H_{f,\CLUST}^{i} \le N$.  If the latter condition holds, then
  $\Sigma$ can be computed from $H_f$ with $Q=\{1,\ldots,\QNUM\}$ and
  with $n_q=\Rank H_{f,\CLUST}^{q}$ for all $q \in Q$.
 \end{Theorem}
  Theorem \ref{theo:pos1} above says 
  that $f$ can be realized by a $K-N$ realization if
  the columns of the Hankel-matrix $H_f$ can be
  divided in $K$ clusters such that for each fixed $l \in I$, all
  columns indexed by $(i,l)$ end up in the same cluster. Moreover, 
  the linear span of the 
  elements of each cluster should have dimension at most $N$.

 In order to prove Theorem \ref{theo:pos1}, we need the
 following transformation of a \LSSJ\ to a $K-N$ realization.
 \begin{Definition}[\LSSJ\ to \PWLS]
\label{lss2pwls}
  Consider 
  a partition $\CLUST=(\CLUST_{l})_{l=1}^{K}$ of $I_f$ for some $K > 0$. 
  Consider an observable \LSSJ\ $\Sigma_l=(n,C,A,\mathcal{X}_0)$ such that
  $\Sigma_l$ is a realization of $f$. Then there exists a unique
  state trajectory $\xi$ of $\Sigma_l$, such that the switching times of
  $\xi$ coincide with the points of non-analyticity $(t_i)_{i=0}^{N_f}$ 
  of $f$ and 
  $\upsilon_{\Sigma_l}(\xi(t))=f(t)$ for all $t \in T$. Define
  the \PWLS\ realization $PW(\Sigma_l,\CLUST)$ of $f$ \emph{associated with 
  $\Sigma_l$ and $\CLUST$} as follows. 
  The \PWLS\ $PW(\Sigma_l,\CLUST)$ is of the form
  \eqref{phas:def}, such that $Q=\{1,\ldots,K\}$, and
  \begin{itemize}
  \item For all $q \in Q$, let  $\hat{\mathcal{X}}_{q,0}$ be the
        set of states $x(t_l)$ such that for some $l \in \CLUST_q$,
        $\xi(t_l)=(1,x(t_l))$. Define
        \[ \mathcal{X}_q=\SPAN\{ A^{k}x \mid k=0,\ldots,n-1,x \in \hat{\mathcal{X}}_{q,0} \} \]
         Set $n_q=\dim \mathcal{X}_q$ and choose a basis of $\mathcal{X}_q$.
  \item
        For each $q \in Q$,  define $A_q$ as the matrix in the basis of
        $\mathcal{X}_q$ of the linear map 
        $\mathcal{X}_q \ni x \mapsto Ax \in \mathcal{X}_q$ obtained
        by restricting $A$ to $\mathcal{X}_q$.
  \item Let $\mathcal{X}_{q,0}$ be the set of vector representations in
        the basis of $\mathcal{X}_q$ of
        the elements of $\hat{\mathcal{X}}_{q,0}$.
  \item Let $C_q$ be the matrix representation in the basis of $\mathcal{X}_q$ 
       of the linear map $\mathcal{X}_q \ni x \mapsto Cx \in \mathbb{R}^p$
       obtained by restricting $C$ to $\mathcal{X}_q$.
  \end{itemize}
 \end{Definition}
 \begin{Proposition}
  If $\Sigma_l$ is a \LSSJ\ realization of $f$, then
  the \PWLS\ $PW(\Sigma_l,\CLUST)$ is a realization of $f$.
 \end{Proposition}
 \begin{pf}[Sketch of the proof of Theorem \ref{theo:pos1}]
  \textbf{only if \ }
    Assume that $\Sigma$ is a \PWLS\ realization of $f$ of
    the form \eqref{phas:def}, such that $|Q| \le K$ and
    $n_q \le N$ for all $q \in Q$. Then there exists a 
    state trajectory $\xi$ of $\Sigma$ such that
    $f=\upsilon_{\Sigma} \circ \xi$ and the switching times of $\xi$ contain
    the points of non-analyticity $(t_i)_{i=0}^{N_f}$ of $f$. Define
    $\CLUST_q=\{i \in I_f \mid \xi(t_i)=(q,x_i) \mbox{ for some } x_i \in \mathcal{X}_{q,0}\}$.
   It then follows that $\CLUST=(\CLUST_q)_{q=1}^{\QNUM}$ is a 
   partitioning of $I_f$. Moreover, for any
   $i \in \CLUST_q$, 
   $M_i^{f}$ is the Markov-parameter of the linear system $(C_q,A_q)$ from
   some initial condition. Hence, by using linear systems theory
   we can show that $\Rank H_{f,\CLUST}^{q} \le n_q \le N$.
  
  \textbf{if\ }
   If the conditions of the theorem hold, then $\Rank H_f = n < +\infty$.
   Construct the \LSSJ\ realization $\Sigma_f$ of $f$, as
  in the proof of Theorem \ref{theo:ex}. Apply Definition
  \ref{lss2pwls} to obtain a linear \PWLS\ $\Sigma=PW(\Sigma_f,\CLUST)$.
   Then $\Sigma$ has $\QNUM \le K$ discrete states.
   Notice that the space $\mathcal{X}_q$ from Definition \ref{theo:pos1} is 
   then isomorphic to the column space of $H_{f,\CLUST}^{q}$,
   $q=1,\ldots,\QNUM$. Hence $n_q = \dim \mathcal{X}_q =\Rank H_{f,\CLUST}^{q} \le N$.
 \end{pf}

Similarly to the general case, in Algorithm \ref{alg05}
 we state an algorithm for computing a 
$K-N$ realization of $f$. 
   \begin{algorithm}
   \caption{
    \newline
    \textbf{Inputs:} Hankel-matrix $H_{f,L,M+1,R}$.  
    \newline
    \textbf{Output:} \PWLS\ $\Sigma$
    }
    \label{alg05} 
    \begin{algorithmic}[1]
    \STATE Apply Algorithm \ref{alg0} to $H_{f,L,M+1,R}$ and 
           denote the result by $\Sigma_{L,M,R}$.
    \STATE 
   \label{alg05:clust}
           Compute a partitioning $\CLUST=(\CLUST_i)_{i=1}^{\QNUM}$  of $I_f^R$
           such that $\QNUM \le K$ and
           for all $i=1,\ldots,\QNUM$, $\Rank H^{i}_{f,\CLUST,L,M,K} \le N$.
           Here $H^{i}_{f,\CLUST,L,M,R}$ is the sub-matrix of $H_{f,L,M,R}$ 
           formed by the columns $(j,l)$ where $j \le M$ and $l \in \CLUST_{q}$. 
    \STATE Return the \PWLS\ $\Sigma=PW(\Sigma_{L,M,R}, \CLUST)$.
    \end{algorithmic}
   \end{algorithm}
\begin{Theorem}[Correctness of Algorithm \ref{alg05}]
  \label{alg05:theo}
   If $f$ has a $K-N$ realization of dimension at most $(K,\min\{L,M\})$, and
   and with a finite set of initial states, then for some $R > 0$,
   Algorithm \ref{alg0} returns a minimal realization of $f$.
  \end{Theorem}

\subsection{Equivalence of \PWLS\ and switched AR models}
\label{sect:sarx}
 Below we show that $f$ is realizable by a \PWLS\ 
 if and only if $f$ satisfies a switched AR model. 
 \begin{Definition}[\SARS\ models]
  A \emph{switched AR system} (abbreviated as \SARS) is a tuple
  \begin{equation} 
  \label{sars:eq1}
    \mathcal{I}=\SAR 
  \end{equation}
  where $Q=\{1,\ldots,\QNUM\}$, $\QNUM > 0$, $n_q > 0$ and
  $A_{q,i} \in \mathbb{R}^{p \times p}$ for all $i=1,\ldots, n_q$.
  The function $f$ is said to satisfy the \SARS\ model
  $\mathcal{I}$, if the following holds. Let $(t_i)_{i=0}^{N_f}$ be
  the points of non-analiticity of $f$. For any $i \in I_f$,
  and for any $t \in [t_i,t_{i+1})$, denote by
  $f^{(k)}(t)$ the $k$th order right-hand derivative
  $\frac{d^k}{ds^k} f(t+s)|_{s=0}$. Then we require that for any $i \in I_{f}$
  there exist $q(i) \in Q$, such that $\forall t \in [t_i,t_{i+1})$,
  \begin{equation}
  \label{sars:eq2}
     f^{(n)}(t) = \sum_{k=1}^{n-1} A_{q(i),k}f^{(n-k)}(t).
  \end{equation}
  The set $Q$ is called the set of discrete modes of $\mathcal{I}$.
 \end{Definition} 
 Notice that if $\QNUM=1$, then the \SARS\ $\mathcal{I}$ corresponds to
 an AR system, and $f$ satisfies $\mathcal{I}$ if it satisfies an AR
 equation. Note, however, that in contrast to trajectories of an AR systems, $f$ is not smooth. 
 \begin{Theorem}
 \label{sars:theo}
  The function $f$ has a realization by a \PWLS\ with $\QNUM$
  discrete states if and only if there exists a \SARS\ $\mathcal{I}$
  with $\QNUM$ discrete modes such that $f$ satisfies $\mathcal{I}$.
 \end{Theorem}
 \begin{pf}[Sketch]
 \textbf{only if\ } 
   Let $\Sigma$ be a realization of $f$. Without loss of
   generality we can assume that
   $\Sigma$ is linear, it is of the form \eqref{phas:def} and
   $n_q=n$ for all $q \in Q$.  Then for each $i \in I_f$, 
   there exists $q(i) \in Q$ such that on $[t_i,t_{i+1})$, $f$
   is an output trajectory of the linear system $(C_{q(i)},A_{q(i)})$.
   From classical theory we then know that there exists 
   matrices $A_{q(i),1},\ldots, A_{q(i),n}$ such that
   \eqref{sars:eq2} holds.
  
 \textbf{if\ }
   Assume $f$ satisfies a \SARS\ $\mathcal{I}$ of the form \eqref{sars:eq2}.
   Define the linear \PWLS\ $\Sigma_{\mathcal{I}}$ of the form \eqref{phas:def}  such that $\Sigma_{\mathcal{I}}$ is a realization of $f$ as follows.
   The sets of discrete modes of $\Sigma_{\mathcal{I}}$ and $\mathcal{I}$
   coincide, $n_q=n$, $\mathcal{X}_{q,0}=\mathbb{R}^{n_q}$ and
   the matrices $C_q$ and $A_{q}$, $q \in Q$ are the matrices of
   the linear system which corresponds to 
   the AR $y^{(n)}(t)=\sum_{j=1}^{n} A_{q,j}y^{(n-j)}(t)$.
   Here we consider the linear system, state vector
   of which is the regressor $x(t)=((y^{(n-1)})^T(t),\ldots, y^T(t))^T$.
   \vspace{-5pt}
 \end{pf}
 Theorem \ref{phas:min} and Theorem \ref{sars:theo} yield
 the following.
 \begin{Corollary}
  The function $f$ has a realization by a \PWLS\ if and only if
  it satisfies an AR system, i.e.
   for some $A_1,\ldots,A_n \in \mathbb{R}^{p \times p}$, $n > 0$, 
   \vspace{-5pt}
  \[ \forall t \in T: f^{(n)}(t) = \sum_{i=1}^{n} A_{i} f^{(n-i)}(t). \]
 \end{Corollary}

\section{Identification algorithm for \PWLS{s}}
\label{sect:ident}
  Below we present an algorithm which
  computes a realization with full observations based on 
  measurements at finitely many time instances. More precisely,
  the algorithm solves the following problem.
 \begin{Problem}
 \label{ident:full_obs}
  Fix integers $\QNUM > 0$ and $n > 0$.
  Consider a piecewise-analytic 
  function $f:T \rightarrow \mathbb{R}^{n}$ and a 
  finite sequence of time instances 
  $t_1 < \ldots < t_M$ and assume that
  $\{f(t_i),\dot f(t_i)\}_{i=1}^{M}$ are known, i.e. we know that
  value of $f$ and its derivative at time instances $t_1,\ldots,t_k$.
  Find a \PWLS\ $\Sigma$ of the form \eqref{phas:def}
  with full-observations such that
  $\Sigma$ realizes $f$, and $C_q=I_n$, $c_q=0$, 
  and $n_q=n$ for all $q \in Q$.
 \end{Problem}
 \textbf{Motivation of the identification problem}
  The motivation for considering the identification problem with
  full observations is the following.
  \begin{enumerate}
  \item
   Notice that the identification problem with full observations and
   the identification problem for \SARS{s} are equivalent. Indeed,
   a \PWLS\ with full observation can be considered as a \SARS{s}
  with $n=2$. Conversely, for a \SARS\ $\mathcal{I}$, the 
  associated \PWLS\ $\Sigma_{\mathcal{I}}$ from the proof of
  Theorem \ref{sars:theo} can be viewed as a \PWLS\ with full
  observations, if the high-order derivatives of $f$ can be measured.
  Hence, by Theorem \ref{sars:theo},
  \emph{the identification problem for
  \PWLS{s} is equivalent to the identification problem for \PWLS\ with
  full observations}, if we assume that high-order derivatives of the
  output can be measured too.
  \item
    The problem of identification with full observations is still a 
    non-trivial problem. Even in this case we have problems
    with identifiability, see Example \ref{example1}.
 \item
   The solution of Problem \ref{ident:full_obs} enables us to
   prove experimentally the feasibility of approximating
   complex systems by \PWLS{s}. This is done by applying the
   identification algorithm which solves Problem \ref{ident:full_obs}
   to simulated trajectories of several well-known complex systems.
  \end{enumerate}
\begin{Example}
\label{example1}
 Consider the \LSSJ{s} $\Sigma_1=(2,C_1,A_1,\mathcal{X}_{0}^1)$ and
 $\Sigma_2=(2,C_2,A_2,\mathcal{X}_0^2)$ where
 $C_1=C_2=I_2$ is the identity matrix and
 $\mathcal{X}_0^{1}=\mathcal{X}_{0}^{2}=\{(1,0)^T\}$ and
 the remaining parameters are as follows;
  \( A_1=\begin{bmatrix} 0 & 0 \\
                         0 & 1 
         \end{bmatrix}
  \) and 
  \( A_{2} = \begin{bmatrix} 0 & 2 \\ 0 & 3 \end{bmatrix} \).
  It then follows that $f(t)=(1,0)^T$, $t \in T$ can be realized
  both by $\Sigma_1$ and $\Sigma_2$, but clearly $\Sigma_1 \ne \Sigma_2$.
  In other words, \emph{realizations of $f$ with full observation are 
  not identifiable.}
  The reason for this phenomenon is that the state of both $\Sigma_1$
  and $\Sigma_2$ live in the subspace $(x,y)$, $y=0$, 
  hence the difference in $A_1$ and $A_2$
  is not visible.
\end{Example}
\begin{Remark}[Related work on \SARS{s}]
  There is a wealth of results on identification of \SARS, mostly
  addressing the discrete-time SISO case. Many of these results, in particular,
  the algebraic approach \cite{MaVidal,VidalAutomatica,Bako09-SYSID2} 
  could probably
  be adapted to the system class of this paper.
  Investigating such adaptations remains a topic of future research.
  Note that it is not at all clear that the convergence results for the discrete-time case have meaningful counterparts in the continuous-time case.
\end{Remark}

In Algorithm \ref{alg2} we present a solution for
solving Problem \ref{ident:full_obs}.
Algorithm \ref{alg2} is an iteration consisting 
of the following steps. 
First, the algorithm starts with a random initialization  of
$\{A_q, a_q\}_{q \in Q}$.
Subsequently, in each iteration step, the current estimates $A_q,a_q, \{w_{q,i}\}_{q \in Q,i=1,\ldots,M}$ are updated by applying first the step 
\ref{alg:PWL-1} for estimating the updated weights with the 
fixed system parameters $\{A_q,a_q\}_{q \in Q}$, 
and then applying step \ref{alg:PWL-2} 
with the previously updated weights $\{w_{q,i}\}_{q \in Q,i=1,\ldots,M}$
 for updating the estimates of the system parameters $\{A_q,a_q\}_{q \in Q}$.
The interpretation of the weights is as follows:
$w_{q,i}$ is one if in $t_i$ the discrete state $q \in Q$
is active, and it is zero otherwise.
The algorithm terminates when an absolute criterion $E$ falls below a pre-specified threshold $\epsilon$. 
The iteration fails if after a pre-defined maximum number of iterations 
$T_{max}$ the criterion $E$ has not yet reached the lower threshold $\epsilon$. 
Unfortunately, the conditions under which the algorithm terminates and
returns a correct realizations are not known yet.
\begin{Remark}[Convergence of Algorithm \ref{alg2}]
  The criterion $E$ of the algorithm always
  decreases or stays the same during the iterations, and
  the same holds for the cost functions in step 
  \ref{alg:PWL-1}--\ref{alg:PWL-2}. 
  Hence, the criterion $E$ and the cost functions in
  step \ref{alg:PWL-1}--\ref{alg:PWL-2} will converge to a fixed value.
  Whether this value is a global optimum remains a topic of future research.
  Note that the true parameters of a realization of $f$ render $E$ zero,
  but there might be many of them, see Example \ref{example1}. 
\end{Remark}

  \begin{algorithm}
  \caption{
   Identification algorithm 
   \newline
   \textbf{Inputs: } data $\{f(t_i),\dot f(t_i)\}_{i=1}^{M}$ 
   \newline
   \textbf{Output: } \PWLS\ $\Sigma$ 
\label{alg2}
  }
  \begin{algorithmic}[1]
  \STATE 
       Initialize the estimates of
       $A_{q},a_{q}$, $q \in  Q$, $k:=0$.
       
   \REPEAT
      \STATE 
      \label{alg:PWL-1}
             Solve the optimization for the weights
             \begin{equation*}
             \begin{split}
              & w = \mbox{arg min}_{w^{*}=\{w^{*}_{q,i}\}_{q \in  Q, i=1,\ldots,M}} E(\{A_{q},a_q\}_{q \in Q}, w^{*}) \\
             & \forall m=1,\ldots,M: \sum_{q \in Q} w^{*}_{q,m} = 1 
               \mbox{, } \forall q \in Q: w^{*}_{q,m} \in [0,1]
             \end{split}
             \end{equation*}
             where for $w=\{w_{q,i}\}_{q \in Q, i=1,\ldots,M}$,
             $A_q \in \mathbb{R}^{n \times n}$, $a_q \in \mathbb{R}^{n}$,
             $q \in Q$, 
		\[
            \begin{split}
            & E(\{A_{q},a_q\}_{q \in Q},w) =\\
            & = \frac{1}{M} \sum_{q \in Q} \sum_{m=1}^{M} w_{q,m} ||\dot f(t_m)  - (A_{q}f(t_m)+a_q)||^{2}
           \end{split}
             \]


             Note that the optimization problem above is a linear
             programming problem, hence the optimal values of
             $w=\{w_{q,i}\}_{q \in Q, i=1,\ldots, M}$ take values
             in the set $\{0,1\}$.
 
      \STATE
       \label{alg:PWL-2}
              Using the weights $w=\{w_{q,i}\}_{q \in Q, i=1,\ldots, M}$
              recalculate the optimal values of $A_{q},a_q$, $q \in Q$ by
              solving the following minimization problem.
              \begin{equation}
              \label{alg:PWL-2:eq1}
           \{A_q,a_q\}_{q \in Q} = \textrm{arg min}_{A^{*}_q,a_q^{*}, q\in Q} E(\{A^{*}_q,a_q^{*}\}_{q \in Q},w)
              \end{equation}
              The optimization problem \eqref{alg:PWL-2:eq1}
              is equivalent to the following linear least squares problem.
              Collect the unknowns $A_{q},a_q$, $q \in Q$ into the matrix $S$
              \begin{eqnarray*}
		  S&=&\begin{bmatrix} A_1 & a_1 & A_2 & a_2 & \cdots & A_{\QNUM} & a_{\QNUM} \end{bmatrix}
              \end{eqnarray*}
              Define the matrices $L \in \mathbb{R}^{(n+1)\QNUM \times M}$ and 
              $Y \in \mathbb{R}^{n \times M}$ 
              \begin{eqnarray*}
                  Y&=&\begin{bmatrix} \dot f(t_1) & \cdots & \dot f(t_m) 
                    \end{bmatrix} \\
                  L &=& \begin{bmatrix}
         w_{11}r_1 & \cdots & w_{M1}r_M \\
         \vdots   & \cdots  & \vdots \\
         w_{1\QNUM}r_1 & \cdots & w_{M\QNUM}r_M
       \end{bmatrix}  \\
               r_i&=&\begin{bmatrix} f(t_i)^T & 1 \end{bmatrix}^T,
               \forall i=1,\ldots,M
              \end{eqnarray*}
             Then the solution of
            \eqref{alg:PWL-2:eq1} can be read off from the entries of $S$,
            where 
              \begin{equation*}
	       S=\mbox{arg min}_{\hat{S} \in \mathbb{R}^{n \times (n+1)\QNUM}} 
               ||Y-\hat{S}L||^{2}
              \end{equation*}

       \STATE 
             Update $E=E(\{A_q,a_q\}_{q \in Q}, w)$ using the new values.

  \UNTIL{$E < \epsilon$}
 \STATE Return $\Sigma$ of the form \eqref{phas:def} with $n_q=n$,
         $C_q=I_n$, $c_q=0$ and
         $A_q,a_q$ as computed in the last iteration, for all $q \in Q$,
         and $\mathcal{X}_{q,0}=\mathbb{R}^{n}$ for all $q \in Q$.
  \end{algorithmic}
  \end{algorithm}

\subsection{Numerical example}
\label{sect:num}
We conducted several numerical experiments.
They were all performed on a
PC with Intel Pentium 4 CPU 2.8 GHz processor and 1.21 GB RAM memory under
 Windows XP, using Matlab 7.5 including the optimization toolbox. In line 
with the definitions above, we use $K$ and $N$ to quantify the
dimension of the \PWLS\ of interest and $M$ to denote the number
of data points.
Dynamical data 
$\{f(t_i),\dot f(t_i)\}_{i=1}^{M}$ was sampled from a given dynamical model 
$\dot x = f(x)$ at regular time intervals $t_i = i.\Delta$, $i=1,\ldots,M$ for
some $\Delta > 0$.  
Several initial states  were chosen randomly, and a 
sampled trajectory was obtained for each of them. We then concatenated
these sampled trajectories 
into one time series. Notice that concatenation of
finite components of two output trajectories of a \PWLS\ is itself a 
finite component of a valid output trajectory. 
Hence, the combined data can be viewed as originating from the 
measurements of one valid output trajectory of a \PWLS.
Gaussian noise with zero mean 
and variance $\sigma^2$ was added to the obtained data.
Finally, Algorithm \ref{alg2} was applied 
to the resulting data. 
In the numerical experiments we first tested the 
performance of the algorithm on artificial \PWLS\ systems. Furthermore, we 
studied the algorithm on Lorenz systems and on the Tyson-Novak model for the cell cycle of budding yeast.

\vspace{-10pt}

\subsubsection*{Simulations on artificial \PWLS\ systems}
We generate an artificial \PWLS\ system containing $K$ affine 
subsystems, each of the same dimension $n_q = N$. The state switching is 
obtained by partitioning the state space in $K$ subsets $\{V_1,\ldots,V_K\}$. 
The switching is generated by using the following switching law:
the discrete mode $i$ is active, if the continuous state belongs to $V_i$.

In the absence of noise, and if the subsystems are sufficiently sampled, the estimation algorithm is usually able to find the system parameters, though occasionally it gets stuck in a local minimum. 
The amount of empirical data required to reconstruct 
the system parameters, within the bounds set by the noise, increases
 proportionally with $N$ and $K$.
Figure \ref{N1} below shows the result of the algorithm for 
$N = 2$, $K = 5$, $M = 1800$ data points. The Voronoi partitioning of state space of the original dataset (black) and the reconstructed dataset (red) are identical. The reconstructed system parameters and the original system parameters
are identical as well.

\begin{figure}
\caption{\label{N1}}
\includegraphics[scale=0.1]{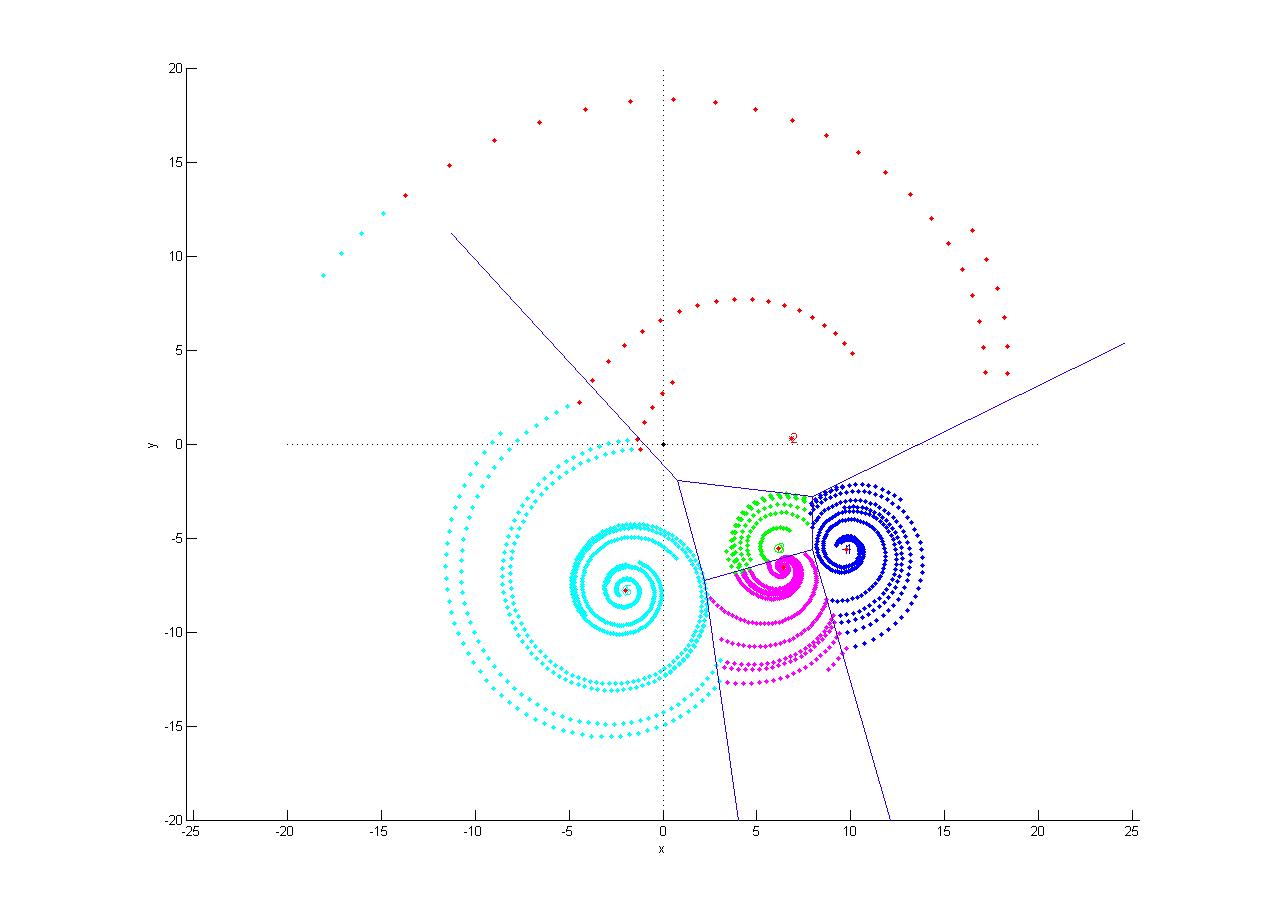}
\end{figure}
Further tests reveal that Algorithm \ref{alg2} is
able to find the correct parameters even if the regions
$\{V_1,\ldots,V_K\}$ are disconnected.



When applying zero-mean Gaussian noise $\mathcal{N}(0,\sigma^2)$,
the accuracy of the reconstruction decreases as $\sigma$ increases.
Figure \ref{N3} below shows the result on a dataset of dimension $N = 2$ with $K = 5$ subsystems and $M = 1800$ data points with a SNR (signal-to-noise ratio) of $5\%$. Two subsystems could not be reconstructed, but the correlation between the three reconstructed systems with the best fitting original subsystems is $98.4$


\begin{figure}
\caption{\label{N3}}
 \includegraphics[scale=0.2]{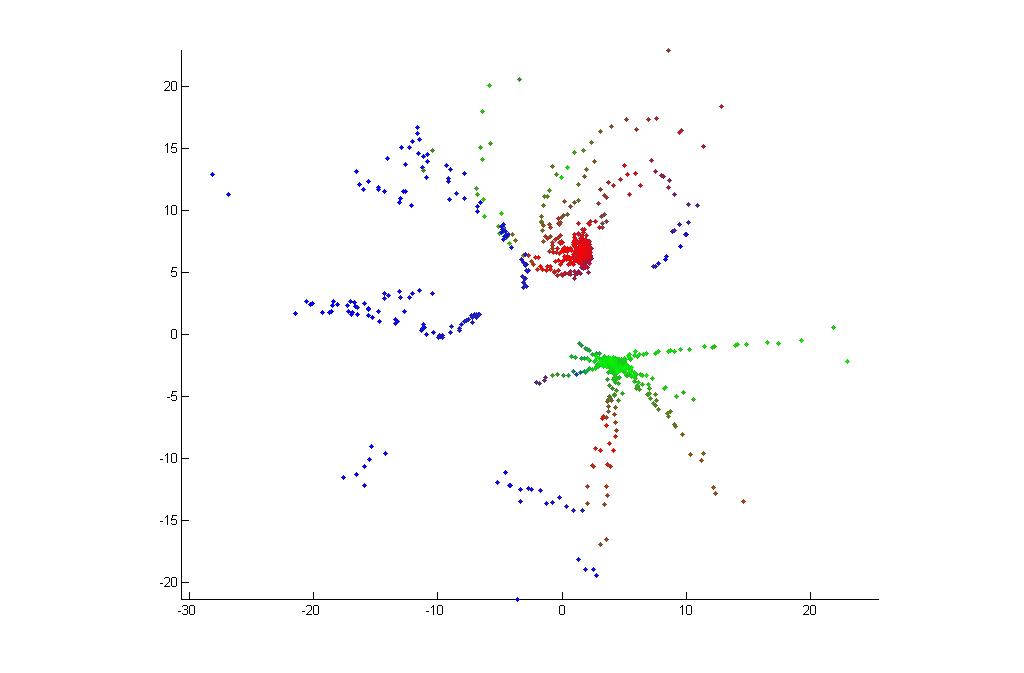}
\end{figure}

\subsubsection{Chaotic systems}

We applied Algorithm \ref{alg2} to Lorenz systems \cite{ref-Lor}.
This example provides insight how the identification 
algorithm performs on systems that are not piecewise affine by nature. 

\begin{figure}
\caption{\label{N4}}
 \includegraphics[scale=0.1]{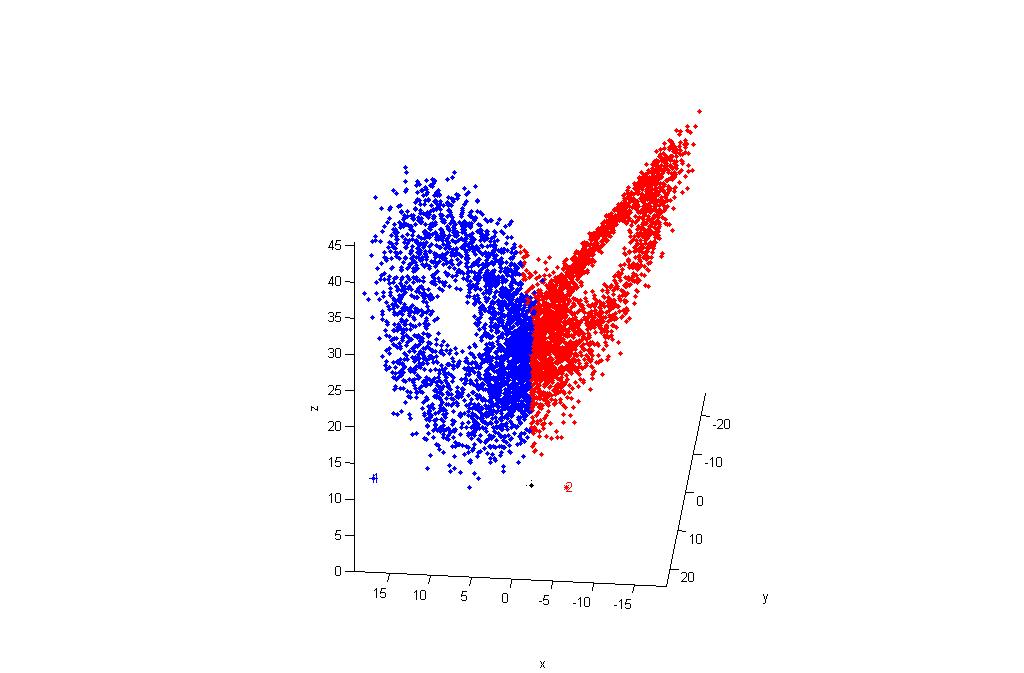}
 \includegraphics[scale=0.1]{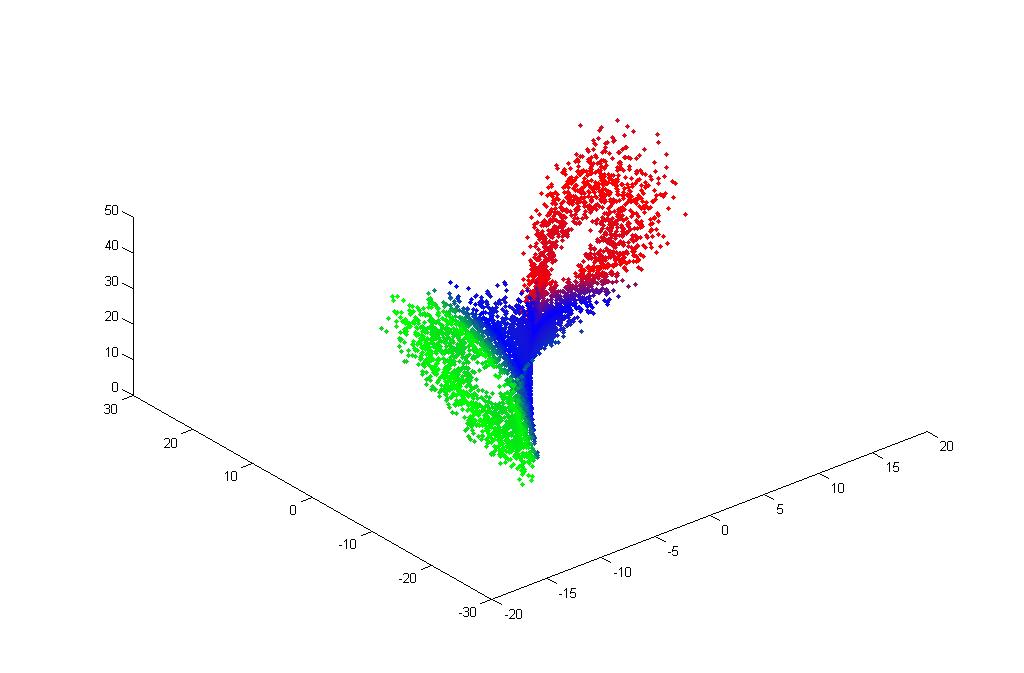}
\end{figure}

The Lorenz attractor can be visualized approximately as consisting of two 
linear subsystems ‘Left’ and ‘Right’, with two different central 
attractors.
A third region can be defined that involves 
the transition between the Left and Right subsystem. In Figure \ref{N4} 
the application of Algorithm \ref{alg2} to Lorenz systems 
with $M=5000$ data points is presented for $K=2,3$.
For $K=2$ the two main parts of the Lorenz attractor are found.  For $K=3$ also the transition region from the left to the right attractor is found. Applying more than $3$ subsystems does not improve the performance.

\subsubsection*{Biological cell cycle models}

\begin{figure}
\caption{\label{N6-a}}
 \includegraphics[scale=0.2]{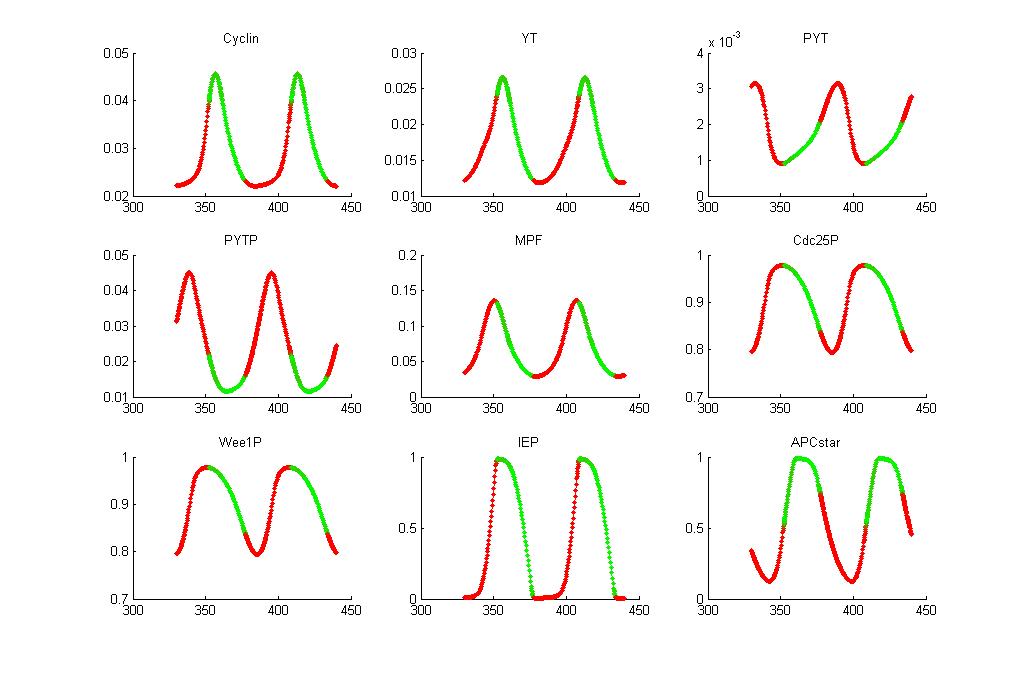}
\end{figure}

We applied Algorithm \ref{alg2} to data simulated using the Tyson-Novak 
model \cite{ref-TN3} for the yeast cell cycle.
 The cell cycle is the succession of events whereby a cell grows and divides into two daughter cells.
The simulated data was used to construct a \PWLS\ approximation
 of the system using Algorithm \ref{alg2}. This resulted in a clear 
partition of the cell cycle dynamics as a function of the selected number 
of classes $K$. The best qualitative results were obtained for $K=2$,
see Figure \ref{N6-a}.
Moreover, for $K > 3$ the results did not improve.
It was found that the obtained two discrete modes correspond to the phases
\emph{S} (DNA synthesis) and \emph{M} (mitosis) of the cell cycle.
The reconstruction remains valid under noise with SNR at most $10\%$. 

\section{ Conclusions }
 We presented some basic results on realization theory of
 \PWLS{s} and a practical identification algorithm. Analysis of
 the correctness and convergence of the presented algorithm remains
 a topic of future research.


\end{document}